# RESEARCH ANNOUNCEMENT

## NOT ALL FREE ARRANGEMENTS ARE $K(\pi, 1)$

PAUL H. EDELMAN AND VICTOR REINER

ABSTRACT. We produce a one-parameter family of hyperplane arrangements that are counterexamples to the conjecture of Saito that the complexified complement of a free arrangement is $K(\pi, 1)$. These arrangements are the restriction of a one-parameter family of arrangements that arose in the study of tilings of certain centrally symmetric octagons. This other family is discussed as well.

### I. DEFINITIONS AND INTRODUCTION

Let $\mathcal{A}$ be a finite set of hyperplanes (subspaces of codimension one) passing through the origin in $\mathbb{R}^d$. The *complexification* of the arrangement $\mathcal{A}$ is the arrangement of hyperplanes in $\mathbb{C}^d$ defined by

$$\mathcal{A}_{\mathbb{C}} = \{ H \otimes_{\mathbb{R}} \mathbb{C} \mid H \in \mathcal{A} \}.$$

Let $\mathcal{M}(\mathcal{A})$ be the complement of $\mathcal{A}_{\mathbb{C}}$ in $\mathbb{C}^d$. We will say that $\mathcal{A}$ is $K(\pi, 1)$ if the space $\mathcal{M}(\mathcal{A})$ is a $K(\pi, 1)$ space; i.e., the universal covering space of $\mathcal{M}(\mathcal{A})$ is contractible and the fundamental group $\pi_1(\mathcal{M}(\mathcal{A})) = \pi$. If $\mathcal{A}$ is $K(\pi, 1)$, then it is known that the cohomology ring $H^*(\mathcal{M}(\mathcal{A}), \mathbb{Z})$ coincides with the group cohomology $H^*(\pi, \mathbb{Z})$.

The *braid arrangement* $\mathcal{A} = A_{d-1}$ is the hyperplane arrangement whose hyperplanes are defined by the linear forms $\{x_i - x_j = 0 \mid 1 \leq i < j \leq d\}$. In 1962 Fadell, Fox, and Neuwirth [FaN, FoN] showed that $\mathcal{A}$ is $K(\pi, 1)$ where $\pi$ is the *pure braid group* on $d$ strands. Subsequently, Arnold [Ar] gave a simple presentation of the cohomology ring, thereby computing the cohomology of the pure braid group. He conjectured that there was a similar presentation of $H^*(\mathcal{M}(\mathcal{A}), \mathbb{Z})$ for an arbitrary arrangement.

Brieskorn [Br] proved this conjecture in 1971 and in 1980 Orlik and Solomon [OS] used these results to give a combinatorial presentation of $H^*(\mathcal{M}(\mathcal{A}), \mathbb{Z})$. Brieskorn also conjectured that all *Coxeter arrangements* are $K(\pi, 1)$. A Coxeter arrangement is the set of reflecting hyperplanes of a finite group acting on $\mathbb{R}^d$ generated by reflections (see, e.g., [Hu]). In particular the braid arrangement is

---

the Coxeter arrangement for the symmetric group $S_d$ permuting the coordinates in $\mathbb{R}^d$. Brieskorn proved the latter conjecture for many Coxeter groups, and it was settled in the affirmative by Deligne [De]. In fact Deligne proved the stronger result that if an arrangement $\mathcal{A}$ is *simplicial*, i.e., every connected component of $\mathbb{R}^d - \cup_{H \in \mathcal{A}} H$ is a union of open rays emanating from the origin whose cross-section is an open simplex, then $\mathcal{A}$ is $K(\pi, 1)$. (We should note that the condition of being simplicial is **not** generic.) In a different direction, Falk and Randell [FR] and Terao [Te] showed that a class of arrangements called *supersolvable arrangements* is also $K(\pi, 1)$.

A common generalization of Coxeter arrangements and supersolvable arrangements is a *free arrangement*, which we now define. For each hyperplane $H$ in $\mathcal{A}$, let $l_H$ be the linear form in the polynomial ring $S = \mathbb{R}[x_1, \ldots, x_d]$ which vanishes on $H$ (so that $l_H$ is uniquely defined up to a scalar multiple). The *module of $\mathcal{A}$-derivations* $D(\mathcal{A})$ is defined to be the set of all derivations $\theta : S \to S$ with the property that $\theta(l_H)$ is divisible by $l_H$ for all $H$ in $\mathcal{A}$. $D(\mathcal{A})$ is a module over the polynomial ring $S$, and we say $\mathcal{A}$ is a *free arrangement* if it is a free module over $S$. If $\mathcal{A}$ is a free arrangement in $\mathbb{R}^d$, then there exists a homogeneous basis $\{\theta_1, \theta_2, \ldots, \theta_d\}$ for $D(\mathcal{A})$ and the degrees of these polynomials (with multiplicities) only depend on $\mathcal{A}$. Call this multiset of degrees the *exponents* of the free arrangement. In the case $\mathcal{A}$ is a Coxeter arrangement, these exponents coincide with the usual definition of exponents of a Coxeter group (see, e.g., [Hu, §3.20]).

In 1975 Saito [Sa] conjectured that if $\mathcal{A}$ is a free arrangement, then it is $K(\pi, 1)$. This conjecture does not completely unify what is known about $K(\pi, 1)$ arrangements because there are simplicial arrangements that are not free. Orlik and Terao [OT, p. 10] remark that this has been one of the two motivating conjectures for most of the recent work on hyperplane arrangements (the other is Terao's conjecture that the freeness of an arrangement is dependent only on the combinatorial properties of the arrangement [OT, Conjecture 4.138]). In the next section we describe a one-parameter family of 3-dimensional arrangements that are not $K(\pi, 1)$, thus disproving Saito's conjecture. This one-parameter family arises as a restriction of a one-parameter family of 4-dimensional arrangements that have some unusual properties as well. We will discuss this other family in the last section.

## II. Counterexamples to Saito's conjecture

Consider the hyperplane arrangement $\mathcal{A}_\alpha$ in $\mathbb{R}^3$ whose hyperplanes are defined by the linear forms

$$\{x, y, z, x-y, x-z, y-z, x-\alpha y, x-\alpha z, y-\alpha z\},$$

where $\alpha \in \mathbb{R}$.

**Theorem 2.1.** *The arrangement $\mathcal{A}_\alpha$ is free for all values of $\alpha$. If $\alpha \neq -1, 0, 1$, then $\mathcal{A}_\alpha$ is not $K(\pi, 1)$.*

*Sketch of proof.* Using the Addition-Deletion Theorem [OT, Theorem 4.51], it is routine to show that $\mathcal{A}_\alpha$ is a free arrangement with exponents

$$\{1, 2, 3\} \text{ if } \alpha = 1 \text{ or } 0,$$
$$\{1, 3, 5\} \text{ if } \alpha = -1,$$
$$\{1, 4, 4\} \text{ otherwise}.$$

Space allowed for art: 14pc

FIGURE 1. THE ARRANGEMENT $\mathcal{A}_{-2}$

To be more explicit, if $\alpha = 1$ or $0$, then $\mathcal{A}_{\alpha}$ is the Coxeter arrangement $A_3$; and if $\alpha = -1$, then $\mathcal{A}_{\alpha}$ is the Coxeter arrangement $B_3$. We thank H. Terao for computing the following basis for $D(\mathcal{A})$ for all other values of $\alpha$:

$$\theta_1 = x\frac{\partial}{\partial x} + y\frac{\partial}{\partial y} + z\frac{\partial}{\partial z},$$

$$\theta_2 = x(x-z)(x+\alpha z)(x+\alpha y)\frac{\partial}{\partial x} + y(y-z)(y+\alpha z)(x+\alpha y)\frac{\partial}{\partial y},$$

$$\theta_3 = x(x-z)(x+\alpha z)(x+(\alpha-1)y-\alpha z)\frac{\partial}{\partial x} + \alpha\, y(y-z)(y+\alpha z)(x-z)\frac{\partial}{\partial y}.$$

To show that $\mathcal{A}_{\alpha}$ is not $K(\pi, 1)$, assume that $\alpha < 0$ and $\alpha \neq -1$. Using elementary linear algebra, one can show that there is a connected component of $\mathbb{R}^3 - \cup_{H \in \mathcal{A}_{\alpha}} H$ bounded by the hyperplanes $\{x - z,\ x - \alpha\, y,\ y - \alpha\, z\}$. Moreover the lines of intersection of pairs of these hyperplanes are not contained in any other hyperplane of the arrangement. Such a configuration is called a *simple triangle*. Figure 1 shows a picture of the arrangement $\mathcal{A}_{-2}$ drawn in the real projective plane with $z = 0$ as the hyperplane at infinity. It follows from results of Hattori [Ha; OT, p. 164] that any 3-dimensional arrangement with a simple triangle is not $K(\pi, 1)$.

By allowing $\alpha$ to range through $\mathbb{C}$, we can construct a *lattice isotopy* in the sense of Randell [Ra; OT, Definition 5.27] between $\mathcal{A}_{\alpha}$ and $\mathcal{A}_{\beta}$ for any two real values of $\alpha$ and $\beta$ neither equal to $0$, $1$, or $-1$. This lattice isotopy implies that $M(\mathcal{A}_{\beta})$ is diffeomorphic to $M(\mathcal{A}_{\alpha})$ [Ra; OT, Theorem 5.28], and hence $\mathcal{A}_{\beta}$ is not $K(\pi, 1)$ for any real value of $\beta$ that is not equal to $0$, $1$, or $-1$.   $\square$

## III. A 4-DIMENSIONAL FAMILY

Given any hyperplane $H$ in $\mathcal{A}$, we define the *restriction arrangement* $\mathcal{A}|_H$ to be the arrangement within the subspace $H$ (thinking of $H$ as $\mathbb{R}^{d-1}$) whose hyperplanes are all of the intersections of hyperplanes of $\mathcal{A}$ with $H$. The 3-dimensional family described in §2 was discovered as the restriction of a 4-dimensional family with its own peculiarities. We will describe that family now.

Let $\mathcal{B}_\alpha$ be the arrangement defined by the forms

$$\{x,\,y,\,z,\,w,\,x-y,\,x-z,\,x-w,\,y-z,\,y-w,$$
$$z-w,\,x-\alpha\,z,\,y-\alpha\,z,\,x-\alpha\,w,\,y-\alpha\,w\}.$$

**Theorem 3.1.** *The arrangement* $\mathcal{B}_\alpha$ *is*

(1) *not free if* $\alpha = -1$;

(2) *a free arrangement with exponents* $\{1,\,2,\,3,\,4\}$ *if* $\alpha = 0,\,1$;

(3) *a free arrangement with exponents* $\{1,\,4,\,4,\,5\}$ *if* $\alpha \neq 0,\,1,\,-1$.

*Sketch of proof.* The proof follows again from a routine application of the Addition-Deletion Theorem [OT, Theorem 4.51]. □

Thus $\mathcal{B}_\alpha$ has the property that in the neighborhood of $\alpha = -1$ a family of free arrangements deforms continuously into one that is not free, and hence the set of free arrangements is not Zariski-closed in the space of all arrangements. By the same token, since non-freeness is a generic condition in the space of *all* hyperplane arrangements, assuming that the coefficients of the linear forms defining the hyperplanes are chosen randomly (see [Zi, Corollary 7.6]), we know that the set of free arrangements is not Zariski-open in the space of all arrangements. However, work of Yuzvinsky [Yu] shows that it **is** a constructible set.

The arrangement $\mathcal{A}_\alpha$ is obtained from $\mathcal{B}_\alpha$ by restricting to any of the last four hyperplanes. The arrangement $\mathcal{B}_\alpha$ arises in the study of tilings of centrally symmetric octagons. Details of this relationship will be discussed in a future paper [ER].

School of Mathematics, University of Minnesota, Minneapolis, Minnesota 55455
*E-mail address,* P. Edelman: `edelman@s1.math.umn.edu`
*E-mail address,* V. Reiner: `reiner@s7.math.umn.edu`